\documentclass[12pt]{amsart}

\usepackage{amssymb}

\usepackage{enumerate}

\usepackage{graphicx}

\makeatletter
\@namedef{subjclassname@2010}{%
  \textup{2010} Mathematics Subject Classification}
\makeatother

\newtheorem{theorem}{ Main Theorem}[section]
\newtheorem{thm}{Theorem}[section]

\theoremstyle{definition}

\numberwithin{equation}{section}

\begin{document}


\baselineskip=17pt


\title [On the sums of many biquadrates  in two different ways]{On the sums of many biquadrates  in two different ways}

\author[F. Izadi]{Farzali Izadi}
\address{Farzali Izadi \\
Department of Mathematics \\ Faculty of Science \\ Urmia University \\ Urmia 165-57153, Iran}
\email{f.izadi@urmia.ac.ir}

\author[M. Baghalaghdam]{Mehdi Baghalaghdam}
\address{Mehdi Baghalaghdam \\
Department of Mathematics\\ Faculty of Science \\ Azarbaijan Shahid Madani University\\Tabriz 53751-71379, Iran}
\email{mehdi.baghalaghdam@azaruniv.edu}

\date{}

\begin{abstract}
The beautiful quartic Diophantine equation $A^4+hB^4=C^4+hD^4$, where $h$ is a fixed arbitrary positive integer, has been studied by some mathematicians for many years. Although Choudhry, Gerardin and Piezas presented solutions of this equation for many values of $h$, the
solutions were not known for arbitrary positive integer values of $h$.
In a separate paper (see the arxiv), the authors completely solved the equation for arbitrary values of $h$, and worked out many examples for different values of $h$, in particular for the values which has not already been given a solution. Our method, give rise to infinitely many solutions and also infinitely many parametric solutions for the equation for arbitrary rational values of $h$.  In the present paper, we use the above solutions as well as a simple idea to show that how  some numbers can be written as the sums of two, three, four, five, or more biquadrates in two different ways. In particular we give examples for the sums of $2$, $3$, $\cdots$, and $10$, biquadrates expressed in two different ways.
\end{abstract}



\maketitle

\section{Introduction}
 The beautiful quartic Diophantine equation $A^4+hB^4=C^4+hD^4$, where $h$ is an arbitrary positive integer, has been studied by some mathematicians for many years. The numerical solutions for $75$ integer values of $h\leq 101$ was given by Choudhry \cite{A.C}. Then these solutions were
first extended by Piezas \cite{T.P} for all positive integer values of $h \leq 101$,
and finally by Tomita \cite{S.T} for all positive integer values of $h< 1000$, except
$h = 967$. The lost solution for $h = 967$, was supplied by Bremner as mentioned by Tomita. Currently, by computer search, the small solutions of this Diophantine equation are known for all positive
integers $h<5000$, and $A, B, C, D<100000$, except for
\\

$h= 1198, 1787, 1987$,\\

$2459, 2572, 2711, 2797, 2971$,\\

$3086, 3193, 3307, 3319, 3334, 3347, 3571, 3622, 3623, 3628, 3644$,\\

$3646, 3742, 3814, 3818, 3851, 3868, 3907, 3943, 3980$,\\

$4003, 4006, 4007, 4051, 4054, 4099, 4231, 4252, 4358, 4406, 4414$,\\

$4418, 4478, 4519, 4574, 4583, 4630, 4643, 4684, 4870, 4955, 4999,$
\\

see \cite{S.T}.
\\

Gerardin and Piezas found solutions of this equation when $h$ is given by polynomials of degrees $5$, and $2$,
respectively, (See \cite{T.P},\cite{S.T}). Also Choudhry presented several new solutions of this
equation when $h$ is given by polynomials of degrees $2$, $3$, and $4$. (See \cite{A.C2})
\\

In a seprate paper, by using the elliptic curves theory, we completely solved the above Diophantine equation for every arbitrary rational values of $h,$ in particular, for every arbitrary positive integer values of $h$, (See \cite{I.B}). In this paper by using the results of the previous paper as well as a simple idea, we easily show that how some numbers can be written  as the sums of two, three, four, five, or more biquadrates in two different ways.
\\

Our main result is the following:
\begin{theorem}
The Diophantine equation $\sum_{i=1}^na_i^4=\sum_{i=1}^nb_i^4$, where $n\geq 2$, is a fixed arbitrary natural number, has infinitely many nontrivial positive integer solutions. Furthermore, we may obtain infinitely many parametric solutions for the aforementioned Diophantine equation. This shows that how some numbers can be written  as the sums of many biquadrates in two different ways.
\end{theorem}

Firstly, we prove  the following theorem, which is our main theorem in the pervious paper, then this theorem easily implies the above main result.
 \section{The Diophantine equation $A^4+hB^4=C^4+hD^4$}
\begin{thm} The Diophantine equation :

 \begin{equation}\label{e8}
A^4+hB^4=C^4+hD^4,
\end{equation}

where $h$ is an arbitrary rational number, has infinitely many integer positive solutions. Furthermore, we may obtain infinitely many parametric solutions for the above Diophantine equation as well.
Moreover, we have $A+C=B+D$.
By letting $h=\frac{v}{u}$, this also solves the equation of the form $uA^4+vB^4=uC^4+vD^4$, for every arbitrary integer values of $u$, and $v$.
\\

 We use some auxiliary variables for transforming the above quartic equation to a cubic elliptic curve of the positive rank in the form
 \\

$Y^2=X^3+FX^2+GX+H$,
\\
where the coefficients $F$, $G$, and $H$, are all functions of $h$.
Since the elliptic curve has positive rank, it has infinitely many rational points which give rise to infinitely many integer solutions for the above equation too.
\end{thm}

Proof. Let: $A=m-q$, $B=m+p$, $C=m+q$, and, $D=m-p$, where all variables are rational numbers. By substituting these variables in the above equation we get
\\

$-8m^3q-8mq^3+8hm^3p+8hmp^3=0$.
\\

Then after some simplifications and clearing the case $m=0$, we get:

\begin{equation}\label{e7}
m^2(hp-q)=-hp^3+q^3.
\end{equation}
\\

We may assume that $hp-q=1$, and $m^2=-hp^3+q^3$. By plugging
$q=hp-1$,

 into the equation \eqref{e7}, and some simplifications, we obtain the elliptic curve:
\\

\begin{equation}\label{e1}
m^2=(h^3-h)p^3-(3h^2)p^2+(3h)p-1.
\end{equation}
\\

Multiplying both sides of this elliptic curve by
 $(h^3-h)^2$, and letting
\begin{equation} \label{e2}
Y=(h^3-h)m,
\end{equation}
and
\begin{equation} \label{e3}
 X=(h^3-h)p,
\end{equation}
\\

we get the new elliptic curve

\begin{equation}\label{e4}
E(h): Y^2=X^3-(3h^2)X^2+(3h(h^3-h))X-(h^3-h)^2.
 \end{equation}
\\

Then if for given $h$, the above elliptic curve has positive rank, by calculating $m$, $p$, $q$, $A$, $B$, $C$, $D$, from the relations \eqref{e2}, \eqref{e3}, $q=hp-1$, $A=m-q$, $B=m+p$, $C=m+q$, $D=m-p$, after some simplifications and canceling the denominators of $A$, $B$, $C$, $D$, we obtain infinitely many integer solutions for the Diophantine equation.
\\

 If the rank of the elliptic curve \eqref{e4} is zero, we may replace $h$ by $ht^4$, for an appropriate arbitrary rational number $t$ such that the rank of the elliptic curve \eqref{e4} becomes positive. Then, we obtain infinitely many integer solutions for the  Diophantine equation
$A^4+(ht^4)B^4=C^4+(ht^4)D^4$. Then by multiplying $t^4$, to the numbers $B^4$, $ D^4$,(written as $A^4+h(tB)^4=C^4+h(tD)^4$), we
get infinitely many positive integer solutions for the main Diophantine equation
$A^4+hB^4=C^4+hD^4$ (see the examples.).

Finally, for getting infinitely parametric solutions, we mention that each point on the elliptic curve can be represented in the form $(\frac{r}{s^2},\frac{t}{s^3})$, where $r$, $s$, $t$ $\in \mathbb{Z}$.\\
So if we put
$nP=(\frac{r_n}{s_n^2},\frac{t_n}{s_n^3})$, where the point $P$ is one of  the elliptic curve generators, we may  obtain a parametric solution for each case of the Diophantine equations by using the new point $P'=nP=(\frac{r_n}{s_n^2},\frac{t_n}{s_n^3})$. Also by using  the new points of infinite order and repeating the above process, we may obtain infinitely many nontrivial parametric solutions for each case of the Diophantine equations. (see \cite{L.W}, page $83$, for more information about the computations of $r_n$, $s_n$, $t_n$.) \\
Now the proof of the above theorem is complete. It is interesting to see that $A+C=B+D$, too.
\\

Remark 1. Note that by putting $h=\frac{v}{u}$, we may solve the Diophantine equation of the form $uA^4+vB^4=uC^4+vD^4$, for every arbitrary integer values of $u$, and $v$.
\\

{\bf Proof of Theorem 1.1}. From the theorem 2.1, we know that  the Diophantine equation\\
  $A^4+hB^4=C^4+hD^4$, where $h$ is an arbitrary fixed rational number, has infinitely many positive integer solutions and we may obtain infinitely many nontrivial parametric solutions for the aforementioned Diophantine equation too.

 Now, in the equation $A^4+hB^4=C^4+hD^4$, let us take $h=\pm h_1^4\pm h_2^4\pm+h_3^4\pm \cdots \pm h_{n-1}^4$, where $h_i$ are arbitrary fixed rational numbers, then we get
 \\

\begin{equation} \label{e111}
A^4+ (\pm h_1^4\pm h_2^4\pm+h_3^4\pm \cdots \pm h_{n-1}^4)B^4=C^4+(\pm h_1^4\pm h_2^4\pm+h_3^4\pm \cdots \pm h_{n-1}^4)D^4.
 \end{equation}
 \\

Now by multiplying $h_i^4$, to the numbers $B^4$, $ D^4$, and by writing the positive terms in the one side and the negative terms in the other side, we get $n$ positive terms of fourth powers in the both sides, and then obtain infinitely many nontrivial solutions and infinitely many nontrivial parametric solutions for the Diophantine equation $\sum_{i=1}^na_i^4=\sum_{i=1}^nb_i^4$. Now the proof of the main theorem is complete.
\\

Remark 2. Surprisingly, we may solve the general Diophantine equation $ \sum_{i=1}^n a_ix_{i} ^4= \sum_{j=1}^na_j y_{j}^4 $, by taking $h=\pm a_1h_1^4\pm \cdots \pm a_mh_m^4$ in the equation $A^4+hB^4=C^4+hD^4$.
\\

Now we are going to work out many examples.
\\

Example 1. Sums of $2$ biquadrates in two different ways ($h=1$):
\\

$A^4+B^4=C^4+D^4$.
\\

$h=16=2^4$.\\

$E(16): Y^2=X^3-768X^2+195840X-16646400$.\\

$rank=1$; Generator: $P=(X,Y)=(340,680)$.\\
\\

Points: $2P=(313,-275)$, $3P=(\frac{995860}{729},\frac{-727724440}{19683})$,
 $4P=(\frac{123577441}{302500},\frac{305200800239}{166375000})$.\\

 $(p,m,q)=(\frac{313}{4080},\frac{-55}{816},\frac{58}{255})$, \\
 $(p',m',q')=(\frac{2929}{8748},\frac{-1070183}{118098},\frac{9529}{2187})$,\\
$(p'',m'',q'')=(\frac{123577441}{1234200000},\frac{305200800239}{678810000000},\frac{46439941}{77137500})$.\\

Solutions:
\\

$1203^4+76^4=653^4+1176^4$,
\\

$1584749^4+2061283^4=555617^4+2219449^4$,
\\

$103470680561^4+746336785578^4=713872281039^4+474466415378^4$.
\\

Example 2. Sums of $3$ biquadrates in two different ways:
\\

$X_1^4+X_2^4+X_3^4=Y_1^4+Y_2^4+Y_3^4$.
\\

$h=\frac{39}{16}=\frac{5^4-1^4}{4^4}$.
\\

$E(\frac{39}{16}): Y^2=X^3-\frac{4563}{256}X^2+\frac{5772195}{65536}X-\frac{2433942225}{16777216}$.
\\

$rank=2$; Generators: $P_1=(X,Y)=(\frac{3289}{256},\frac{3289}{256})$, and  $P_2=(X',Y')=(\frac{6565}{256},\frac{43615}{512})$.
\\

Points: $2P_1=(\frac{4069}{256},\frac{-14183}{512})$, $3P_1=(\frac{9572761}{57600},\frac{1752134549}{864000})$, $P_2$.
\\

 $(p,m,q)=(\frac{5008}{3795},\frac{-8728}{3795},\frac{2804}{1265})$,\\
$(p',m',q')=(\frac{605392}{43875},\frac{110806928}{658125},\frac{36712}{1125})$, \\
 $(p'',m'',q'')=(\frac{1616}{759},\frac{488}{69},\frac{1060}{253})$.\\
\\

Solutions:
\\

$8570^4+2325^4+1717^4=158^4+8585^4+465^4$,
\\

$11166301^4+18732470^4+3178939^4=16535431^4+15894695^4+3746494^4$,
\\

$1094^4+4365^4+469^4=4274^4+2345^4+873^4$.
\\

Example 3. Sums of $4$ biquadrates in two different ways:
\\

$X_1^4+X_2^4+X_3^4+X_4^4=Y_1^4+Y_2^4+Y_3^4+Y_4^4$.
\\

$h=23=\frac{5^4-1^4-4^4}{2^4}$.\\

$E(23): Y^2=X^3-1587X^2+837936X-147476736$.\\

rank=1; Generator: $P=(X,Y)=(880,6512)$.\\

Points: P, $2P=(\frac{3424933}{5476},\frac{275924489}{405224})$.\\

$(p,m,q)=(\frac{5}{69},\frac{37}{69},\frac{46}{69})$,\\
$(p',m',q')=(\frac{3424933}{66500544},\frac{275924489}{4921040256},\frac{533605}{2891328})$.\\
 \\

Solutions:
\\

$9^4+105^4+16^4+64^4=83^4+80^4+21^4+84^4$,
\\

$1264542442^4+2646847655^4+22479447^4+89917788^4=\\
2368240398^4+112397235^4+529369531^4+2117478124^4$.
\\

Example 4. Sums of $5$ biquadrates in two different ways:
\\

$X_1^4+X_2^4+X_3^4+X_4^4+X_5^4=Y_1^4+Y_2^4+Y_3^4+Y_4^4+Y_5^4$.
\\

 $h=\frac{3}{17}=\frac{1^4+2^4+3^4+11^4}{17^4}$.\\

$E(\frac{3}{17}): Y^2=X^3-\frac{27}{289}X^2-\frac{7560}{83521}X-\frac{705600}{24137569}$.\\

rank=1; Generator: $P=(X,Y)=(\frac{400}{289},\frac{440}{289})$.\\

Points: P, $2P=(\frac{65437}{139876},\frac{313237}{3077272})$.\\

 $(p,m,q)=(\frac{-170}{21},\frac{-187}{21},\frac{-51}{21})$,\\
$(p',m',q')=(\frac{-1112429}{406560},\frac{-5325029}{8944320},\frac{-200957}{135520})$. \\
 \\

Solutions:

 $2312^4+357^4+714^4+1071^4+3927^4=4046^4+17^4+34^4+51^4+187^4$,
 \\

$134948261^4+29798467^4+59596934^4+89395401^4+327783137^4=\\
315999247^4+19148409^4+38296818^4+57445227^4+210632499^4$.
\\

Example 5. Sums of $6$ biquadrates in two different ways:
\\

$X_1^4+X_2^4+X_3^4+X_4^4+X_5^4+X_6^4=Y_1^4+Y_2^4+Y_3^4+Y_4^4+Y_5^4+Y_6^4$.
\\

$h=\frac{66}{25}=\frac{1^4+2^4+3^4+4^4+6^4}{5^4}$.\\

$E(\frac{66}{25}): Y^2=X^3-\frac{13068}{625}X^2+\frac{48756708}{390625}X-\frac{60637092516}{244140625}$.\\

rank=1; Generator: $P=(X,Y)=(\frac{1020552759889}{78568090000},\frac{5339057694122399}{880905425080000})$.\\

Point: P.\\

$(p,m,q)=(\frac{47868328325}{58077532128},\frac{250424844940075}{651165290219136},\frac{1034770525}{879962608})$.\\

Solution:
 \\

 $103059413079145^4+31484981684799^4+62969963369598^4+94454945054397^4+
 125939926739196^4+188909890108794^4=203229351055175^4+11450994089593^4+22901988179186^4+
 34352982268779^4+45803976358372^4+68705964537558^4$.
 \\

Example 6. Sums of $7$ biquadrates in two different ways:
\\

$X_1^4+X_2^4+X_3^4+X_4^4+X_5^4+X_6^4+X_7^4=Y_1^4+Y_2^4+Y_3^4+Y_4^4+Y_5^4+Y_6^4+Y_7^4$.
\\

 $h=\frac{77}{3}=\frac{6^4-1^4-2^4-3^4+4^4-5^4}{3^4}$.\\

$E(\frac{77}{3}): Y^2=X^3-\frac{5929}{3}X^2+\frac{35099680}{27}X-\frac{207790105600}{729}$.\\

rank=1; Generator: $P=(X,Y)=(\frac{92500}{81},\frac{7695260}{729})$.\\

Points: P, $2P=(\frac{6892959356452}{8759275281},\frac{975806887820176684}{819789332824071})$.\\

$(p,m,q)=(\frac{125}{1848},\frac{10399}{16632},\frac{53}{72})$,\\
$(p',m',q')=(\frac{1723239839113}{36970630037880},\frac{243951721955044171}{3460118235875227080},\frac{282825681793}{1440414157320})$. \\

 Solutions:
 \\

$2766^4+34572^4+23048^4+28810^4+4637^4+9274^4+13911^4=\\
 33963^4+27822^4+18548^4+23185^4+5762^4+11524^4+17286^4$,
\\

$X_1=653165044877947269$, $X_2=1215694385212406862$, $X_3=810462923474937908$,\\
$X_4=1013078654343672385$, $X_5=41335991086309694$, $X_6=82671982172619388$,\\
$X_7=124007973258929082$, $Y_1=1385020210743079782$, $Y_2=248015946517858164$,\\
$Y_3=165343964345238776$, $Y_4=206679955431548470$, $Y_5=202615730868734477$,\\
$Y_6= 405231461737468954$, $Y_7=607847192606203431$.\\
\\

Example 7. Sums of $8$ biquadrates in two different ways:
\\

$X_1^4+X_2^4+ \cdots +X_8^4=Y_1^4+Y_2^4+ \cdots +Y_8^4$.
\\

$h=10=\frac{7^4+1^4+2^4-3^4-4^4-5^4-6^4}{2^4}$.\\

$E(10): Y^2=X^3-300X^2+29700X-980100$.\\

$rank=1$; Generator: $P=(X,Y)=(165, 495)$.\\

Points: $P$, $2P=(\frac{505}{4},\frac{-85}{8})$, $3P=(\frac{172029}{961}, \frac{-20192733}{29791})$.\\

$(p,m,q)=(\frac{1}{6},\frac{1}{2},\frac{2}{3})$, $(p',m',q')=(\frac{101}{792},\frac{-17}{1584},\frac{109}{396})$, $(p'',m'',q'')=(\frac{5213}{28830},\frac{-203967}{297910},\frac{2330}{2883})$.\\

Solutions:\\
$3^4+4^4+4^4+5^4+14^4=7^4+7^4+8^4+10^4+12^4$,\\

$906^4+1295^4+185^4+370^4+657^4+876^4+1095^4+1314^4=\\
838^4+1533^4+219^4+438^4+555^4+740^4+925^4+1110^4$,\\

$1334201^4+1576043^4+225149^4+450298^4+1160256^4+1547008^4+1933760^4+2320512^4=\\
110399^4+2707264^4+386752^4+773504^4+675447^4+900596^4+1125745^4+1350894^4$.
\\

Example 8. Sums of $9$ biquadrates in two different ways:
\\

$X_1^4+X_2^4+ \cdots +X_9^4=Y_1^4+Y_2^4+ \cdots +Y_9^4$.
\\

$h=\frac{21}{8}=\frac{8^4+1^4-2^4-3^4-4^4+5^4-6^4-7^4}{4^4}$.\\

$E(\frac{21}{8}): Y^2=X^3-\frac{1323}{64}X^2+\frac{498771}{4096}X-\frac{62678889}{262144}$.\\

$rank=1$; Generator: $P=(X,Y)=(\frac{163241}{11552}, \frac{46525193}{3511808})$.\\

Point: $P$.

$(p,m,q)=(\frac{6928}{7581},\frac{123409}{144039},\frac{505}{361})$.\\

Solution:\\

$312344^4+2040328^4+255041^4+1275205^4+16446^4+24669^4+32892^4+49338^4+57561^4=\\
1299616^4+65784^4+8223^4+41115^4+510082^4+765123^4+1020164^4+1530246^4+1785287^4$.
\\

By choosing the other points on the elliptic curve such as $2P$, $3P$, $\cdots$, (or changing the value of $h$, and getting new elliptic curve) we obtain infinitely many solutions for the above  Diophantine equation as well.
\\

$h=\frac{-3}{2}=\frac{8^4+1^4+2^4+3^4-4^4-5^4-6^4-7^4}{4^4}$.\\

$E(\frac{-3}{2}): Y^2=X^3-\frac{27}{4}X^2+\frac{135}{16}X-\frac{225}{64}$.\\

$rank=1$; Generator: $P=(X,Y)=(\frac{85}{16}, \frac{55}{64})$.\\

Point: $P$.

$(p,m,q)=(\frac{-17}{6},\frac{-11}{24},\frac{13}{4})$.\\

Solution:\\

$356^4+632^4+79^4+158^4+237^4+228^4+285^4+342^4+399^4=\\
268^4+456^4+57^4+114^4+171^4+316^4+395^4+474^4+553^4$.\\
\\

Example 9. Sums of $10$ biquadrates in two different ways:
\\

$X_1^4+X_2^4+ \cdots +X_{10}^4=Y_1^4+Y_2^4+ \cdots +Y_{10}^4$.
\\

$h=-63=\frac{14^4+1^4+2^4+3^4+4^4+5^4-6^4-11^4-13^4}{3^4}$.\\

$E(-63): Y^2=X^3-11907X^2+47246976X-62492000256$.\\

$rank=1$; Generator: $P=(X,Y)=(4960, 30752)$.\\

Points: $P$, $2P=(\frac{4096948}{961},\frac{74223316}{29791})$.
\\

$(p,m,q)=(\frac{-5}{252},\frac{-31}{252},\frac{1}{4})$, $(p',m',q')=(\frac{-1024237}{60058656}, \frac{-18555829}{1861818336},\frac{70925}{953312})$.\\

Solutions:
\\

$141^4+252^4+18^4+36^4+54^4+72^4+90^4+78^4+143^4+169^4=\\
48^4+182^4+13^4+26^4+39^4+52^4+65^4+108^4+198^4+234^4$,
\\

$235608531^4+352150232^4+25153588^4+50307176^4+75460764^4+100614352^4+125767940^4+
39586554^4+72575349^4+85770867^4=\\
179941044^4+92368626^4+6597759^4+13195518^4+19793277^4+26391036^4+32988795^4+
150921528^4+276689468^4+326996644^4$.
\\

By choosing the other points on the above elliptic curves such as $3P$, $4P$, $\cdots$, (or changing the value of $h$, and getting new elliptic curves) we obtain infinitely many solutions for the above Diophantine equation.
\\

The Sage software was used for calculating the rank of the elliptic curves, (see \cite{S.A}).


\begin{thebibliography}{HD}
\bibitem[A.C]{A.C} A. Choudhry, : On The Diophantine Equation $A^4+hB^4=C^4+hD^4$, Indian J. Pure Appl. Math,  $26(11), pp. 1057-1061, (1995)$.
\bibitem [A.C2]{A.C2}  A. Choudhry, A Note on the Quartic Diophantine Equation $A^4+hB^4=C^4+hD^4$,
 available at arXiv. ($2016$)

\bibitem[I.B]{I.B}  F. Izadi, and M. Baghalaghdam, ''Is the quaratic Diophantine equation  $A^4 + hB^4 = C^4 + hD^4$, solvable for any integer $h$?'', submitted, ($2017$).
\bibitem[T.P]{T.P} T. Piezas, A collection of algebraic identities, available at
https://sites.google.com/site/tpiezas/0021e, (accessed on 7 April ($2016$)).
\bibitem[S.T]{S.T} S. Tomita, https://www.maroon.dti.ne.jp/fermat/dioph121e.html, ( accessed on 7 April $(2016)$).

\bibitem[S.A]{S.A} Sage software, available from http://sagemath.org.
\bibitem[L.W]{L.W}  L. C.  Washington, Elliptic Curves: Number Theory and Cryptography, Chapman-Hall, $(2008)$



\normalsize
\baselineskip=17pt



\end{thebibliography}
\end{document}